\documentclass{article}

\usepackage{arxiv}

\usepackage[utf8]{inputenc} % allow utf-8 input
\usepackage[T1]{fontenc}    % use 8-bit T1 fonts
\usepackage{hyperref}       % hyperlinks
\usepackage{url}            % simple URL typesetting
\usepackage{booktabs}       % professional-quality tables
\usepackage{amsfonts}       % blackboard math symbols
\usepackage{nicefrac}       % compact symbols for 1/2, etc.
\usepackage{microtype}      % microtypography

\usepackage{lipsum}         % Can be removed after putting your text content
\usepackage{graphicx}
\usepackage[square,numbers]{natbib}
\usepackage{doi}

\usepackage{amsmath}
\usepackage{amssymb}
\usepackage{amsthm}
\usepackage{amssymb}

\usepackage{mathrsfs}

\usepackage{cleveref}       % smart cross-referencing

\DeclareMathOperator{\Q}{\mathbb Q}
\DeclareMathOperator{\R}{\mathbb R}

\DeclareMathOperator{\supp}{Supp}
\DeclareMathOperator{\vcsp}{VCSP}

\DeclareMathOperator{\wpol}{wPol}

\DeclareMathOperator{\g}{\Gamma}

\DeclareMathOperator{\ar}{ar}
\DeclareMathOperator{\QQ}{\mathbb Q\cup\{+\infty\}}
\DeclareMathOperator{\RR}{\mathbb R\cup\{+\infty\}}

\DeclareMathOperator{\imp}{Imp}

\DeclareMathOperator{\lcm}{lcm}

\DeclareMathOperator{\cone}{Cone}

\newcommand{\defeq}{\mathrel{\mathop:}=}

\makeatletter
\newcommand*{\rom}[1]{\expandafter\@slowromancap\romannumeral #1@}
\makeatother

\theoremstyle{plain}
\newtheorem{theorem}{Theorem}[section]
\newtheorem{lemma}[theorem]{Lemma}
\newtheorem{proposition}[theorem]{Proposition}
\newtheorem{corollary}[theorem]{Corollary}

\theoremstyle{definition}
\newtheorem{definition}[theorem]{Definition}

\newtheorem{remark}[theorem]{Remark}

\title{An Application of Farkas' Lemma to\\ Finite-Valued Constraint Satisfaction Problems\\ over Infinite Domains}

\date{}

\author{ \href{https://orcid.org/0000-0002-4126-2758}{\includegraphics[scale=0.06]{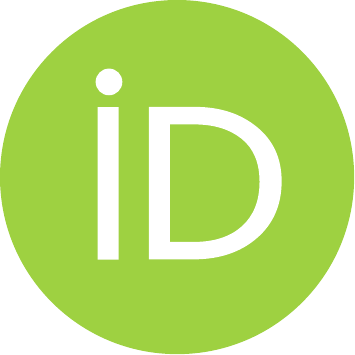}\hspace{1mm}Friedrich Martin Schneider}\\Institute of Discrete Mathematics and Algebra\\ Technische Universit{\"a}t Bergakademie Freiberg\\ Germany\\
	\texttt{martin.schneider@math.tu-freiberg.de} \\
	\And
	\href{https://orcid.org/0000-0002-7312-5002}{\includegraphics[scale=0.06]{orcid.pdf}\hspace{1mm}Caterina Viola}\\Faculty of Mathematics and Physics\\ Charles University Prague\\ Czech Republic\\
	\texttt{caterina.viola@matfyz.cuni.cz}}

\renewcommand{\headeright}{ }
\renewcommand{\undertitle}{ }
\renewcommand{\shorttitle}{An Application of Farkas' Lemma to Infinite-Domain Finite-Valued CSPs}

\begin{document}
	\maketitle
	
	\begin{abstract}
		We show a universal algebraic local characterisation of the expressive power of  finite-valued languages with domains of arbitrary cardinality and containing arbitrary many cost functions.~\footnote{The paper is based on a chapter from the second author's doctoral dissertation~\cite{caterinathesis}.}
	\end{abstract}


%\tableofcontents

\section{Introduction}

Constraint satisfaction problems (CSPs) are computational decision problems in which the input is specified by a bunch of constraints involving some variables and the task is to decide whether it is possible to compute values for the variables that satisfy all given constraints. The constraint satisfaction framework is very expressive and the class of computational problems that can be modelled as CSPs is broad and contains, e.g., the famous \textsc{SAT} and \textsc{3-Colorability}~\cite{Karp1972}. In fact, if we allow the domain to be an infinite set, then every computational decision problem can be reduced\footnote{by a Turing reduction.} to a CSP~\cite{BodirskyGrohe}.

However, there are many situations in which there are ways of satisfying a constraint that are preferable with respect to others, or in which there are no feasibility constraints to satisfy but there are combinations of values for the variables that are more desirable (or less expensive) than others.  \emph{Valued constraint satisfaction problems (VCSPs)} are a generalisation of CSPs that capture  these scenarios. In a VCSP, the constraints are \textit{valued}, that is every possible evaluation of a given tuple of variables comes with a measure of its cost (or desirability), namely a value.
	
Formally, the input of a VCSP is a finite sum of $\RR$-valued cost functions (the objective function) defined on a fixed set called the domain and applied to finitely many variables.  In the case in which the domain is finite, the task is to compute an assignment that minimises the objective function; in the arbitrary-domain setting, a threshold from $\RR$ is given as part of the instance and the task is to compute an assignment of values from the domain for the variables whose cost, i.e., the value of the objective function corresponding to the assignment, is \emph{strictly smaller} than the given threshold. 

By choosing a suitable set of cost functions (the \emph{valued constraint language}), many optimisation problems arising from several scientific fields can be modelled as VCSPs. The list of problems that can be modelled as VCSPs includes  \textsc{Linear Programming}~\cite{Khachiyan,Karmarkar1984}, \textsc{Min Cost Homomorphism}~\cite{Takhanov10:stacs}, \textsc{Min Feedback Arc Set}~\cite{Karp1972}, \textsc{Min Correlation Clustering}~\cite{CorrelationClustering}, and several cutting and clustering problems over graphs~(see, e.g., \cite[Section~2]{krokhin2017complexity}). In turns, these families of problems find applications in areas of artificial intelligence such as, for example, local clustering, community detection \cite{pmlr-v28-buhler13}, segmentation \cite{Klodt-et-al-13}, and 3D reconstruction \cite{Kolev_et_al_pami12}.

By restricting the class of VCSPs to instances such that the cost 
functions take \emph{real-only-values}, i.e., values in $\R$ rather than in $\RR$, we obtain the class of \emph{finite-valued CSPs}. Finite-valued CSPs are pure computational optimisation problems. 
Sometimes, the terminology \emph{general-valued CSPs} is used to distinguish VCSPs with $\RR$-valued cost functions from finite-valued CSPs. Also the corresponding valued languages are referred to respectively as \emph{finite-valued languages} and \emph{general-valued languages}. 

In the case in which the domain of the cost functions is a finite set, the computational complexity of finite-valued CSPs~\cite{ThapperZivnyfinitevalued}, general-valued  CSPs~\cite{GenVCSP15},  and classical CSPs~\cite{BulatovFVConjecture,ZhukFVConjecture} has been completely classified.  To obtain this classification results in the finite-domain setting,  the use of the so-called universal algebraic approach~\cite{JBK,KozikOchremiak15}  was 
decisive. It turned the study of the complexity of (the VCSP for) a general-valued constraint language into the study of certain algebraic objects called \emph{weighted polymorphisms}.

 The present article is the very first attempt to relate the set of weighted polymorphisms of a finite-valued language \textit{over an infinite domain} to the set of cost functions that it can express. This is the first step towards understanding  whether  the computational complexity of VCSPs for finite-valued languages over infinite domains can be classified -- similarly to the finite-domain case -- in terms of the set of weighted polymorphisms of the finite-valued  language, yielding a \textit{topological approach} to the study of the complexity of infinite-domains VCSPs. 
 
 Our contribution is a characterisation of the set of cost functions improved by all weighted polymorphisms improving the cost functions from a given finite-valued language over an  \emph{infinite domain}. We will show that for domains of arbitrary cardinality this set coincides with a local version of the expressive power. The result relies on an application of an infinite-dimensional version of Farkas' lemma proved by Swartz~\cite{swartz}.

\section{Preliminaries and Notation}
Throughout this paper, we %denote by $x_i$ the $i$-th component of a tuple $x$. We
denote by $\R$ ($\R_{\geq0}$, $\R_{>0}$, resp.) the set of the real numbers (non-negative real numbers, positive real numbers, resp.).  We denote by $\Q$ the set of the rational numbers.  
In the remainder, for a map $f \colon A^r \to B$, we define the \emph{arity} of $f$ to be the number $\ar(f)=r$.

Let $M$ be a set. As usual, the \emph{support} of a  function $\alpha \colon M \to \R$ is defined as the set $\supp(\alpha)\defeq \{x \in M \mid \alpha(x)\neq0\}$. We consider the vector space \[\mathbb{R}[M] 
\defeq \left\{ \alpha \in \mathbb{R}^{M} \, \mid \supp (\alpha) \text{ is 
	finite} \right\}.\] For any $x \in M$, let us define $\delta_{x} \in \mathbb{R}[M]$ by \[
\delta_{x}(y) \, \defeq \, \begin{cases}
\, 1 & \text{if } y=x, \\
\, 0 & \text{otherwise}
\end{cases} \qquad (y \in M).
\] Furthermore, let us consider  $\mathbb{R}_{\geq 0}[M] \defeq \{ \alpha \in \mathbb{R}[M] \mid \forall x \in M \colon \, \alpha (x) \geq 
0 \}$.  

Let $D$ be a set. For every $r \in \mathbb{N}$ we define $\mathrm{wRel}_{D}^{(r)} \defeq \mathbb{R}^{D^{r}}$, i.e., $\mathrm{wRel}_{D}^{(r)}$ is the set of all the functions from $D^r$ to $\R$. We also set \[\mathrm{wRel}_{D} \defeq \bigcup\nolimits_{r \in \mathbb{N}} \mathrm{wRel}_{D}^{(r)}.\] Moreover, let $k \in \mathbb{N}\setminus \{ 0 \}$. Then \begin{itemize}
	\item $\mathcal O_D^{(k)} \defeq D^{D^{k}}$ denotes the set of all $k$-ary operations over $D$;
	\item for each $i \in \{ 1, \ldots,k\}$, the \emph{$i$-th projection of arity $k$} is defined as \begin{displaymath}
			\qquad e_{i}^{(k)} \colon \, D^{k} \, \longrightarrow \, D, \quad (x_{1},\ldots,x_{k}) \, \longmapsto \, x_{i} ;
		\end{displaymath}
	\item $\mathcal J_D^{(k)} \defeq \left.\! \left\{ e_{i}^{(k)} \right\vert i \in \{1,\ldots,k\} \right\}\!$ denotes the subset  of $\mathcal O_D^{(k)}$  consisting of all $k$-ary projections.
\end{itemize}

\subsection{Valued Constraint Satisfaction Problems}
\label{sect:vcsps}
Let $D$ be a set. A \emph{cost function} over $D$ is a function $f \in \mathrm{wRel}_{D}^{(\ar(f))}$.
A \emph{finite-valued language $\Gamma$ with domain
	$D$}  is a set of cost functions over $D$.
%%	\item a signature $\tau$ 
%	consisting of function symbols $f$, each equipped with an arity $\ar(f)$, 
%	\item a set $D = \dom(\Gamma)$ (the \emph{domain}), 
%	\item for each $f \in \tau$ a 
%%\end{itemize}

Let $\Gamma$ be a finite-valued language with domain $D$. 
The \emph{valued constraint satisfaction problem  for $\Gamma$}, denoted by $\vcsp(\g)$, is the following computational problem.

\begin{definition}\label{vcspdef}
	An \emph{instance}  of $\vcsp(\g)$ is a triple $I:=(V_I, \phi_I,u_I)$
	where 
	\begin{itemize}
		\item $V_I$ is a finite set of variables, 
		\item $\phi_I$ is an expression  of the form
		\[\sum_{i=1}^{m} f_i(v^i_1,\ldots, v^i_{\ar(f_i)})\]
		where $f_1,\dots,f_m \in \g$ and all the $v^i_j$ are variables from $V_I$, and
		\item $u_I$ is a value from $\Q$. 
	\end{itemize} 
	The task is to decide whether there exists an assignment $\alpha \colon V_I \to D$, whose \emph{cost}, defined as
	\[\phi_I(\alpha(v_1),\ldots, \alpha(v_{\lvert V_I \rvert})):=\sum_{i=1}^{m} f_i(\alpha(v^i_1),\ldots, \alpha(v^i_{\ar(f_i)}))\] 
	is {\it strictly smaller} than $u_I$.
\end{definition}

We remark that we need the value $u_I$ from Definition~\ref{vcspdef} to be in $\Q$ (rather than in $\R$) because we need to represent it computationally.

%\begin{remark}The problem $\vcsp(\g)$ can also be defined for the case in 	which the cost functions in $\g$ are partially defined. In this case, every cost function $f \in \g$ is  a cost function $f\colon D^{\ar(f)}\to \RR$, where $+\infty$ is an extra element with the expected properties that 	for all $c \in {\mathbb Q} \cup \{+\infty\}$	\begin{align*}	(+\infty) + c & = c + (+\infty) = +\infty, \\ 	\text{  and } c & < +\infty \text{ iff } c \in {\mathbb Q},	\end{align*}   and for $x \in D^{\ar(f)}$ having $f(x)=+\infty$ is interpreted as $f$ not defined on $x$, that is, $x \notin \dom(f)$.\end{remark}

Let $\g$ be a general-valued language with domain $D$. For every $C \subseteq D$, 
the general-valued language $\g'$ with domain $C$ whose cost functions are exactly the restrictions to $C$ of the cost functions from $\g$ is called a \emph{general-valued sublanguage} of $\g$.  

We say that a general-valued language has \emph{finite size} if it contains only finitely many cost functions; otherwise, we say that it has \emph{infinite size}.

Note that, given a general-valued language $\g$, if it has finite size then it is 
inessential for the computational complexity
of $\vcsp(\g)$ how the cost functions appearing in $\phi_I$ are represented.
In the case of general-valued languages with infinite size we have to  fix an explicit representation of the cost function, as the computational complexity of the corresponding VCSP depends on such a representation.

Let $\g$ be a finite-valued language with domain $D$ and infinite size. A \emph{finite-valued finite reduct} of $\g$ is a valued language $\g'$ with domain $D$ only containing a finite subset of the cost functions from $\g$.  

%\subsection{Expressive Power and Weighted Polymorphisms for Finite-Domain Finite-Valued Languages}

%The universal algebraic approach for finite-domain VCSPs showed that to obtain a complexity classification it is enough to study just those finite-valued constraint languages that are characterised by having a given set of weighted polymorphisms. More in detail, in the finite-domain setting, the set of cost functions \emph{improved} by all \emph{weighted polymorphisms} improving the cost functions from a given finite-valued language $\Gamma$ coincides with the set of cost functions which are \emph{expressible} in $\Gamma$, that is, the set of cost functions which can be written as the pointwise minimum of a conic combination of cost functions from $\Gamma$. The set of all cost functions expressible in a finite-valued language $\g$ is called the \emph{expressive power} of $\g$. If the domain of the cost functions in $\Gamma$ is a finite set, then the VCSP for the expressive power of $\Gamma$ is polynomial-time many-one reducible to the VCSP for $\Gamma$. 

\subsection{Expressive Power}

\begin{definition}\label{def:exprpwr1}Let $\g$ be a finite-valued language with an arbitrary domain $D$ and let $\rho \in \mathrm{wRel}_{D}$. We say that $\rho$ is an \emph{affine conical combination} over $\g$ if there 
	exist \begin{itemize}\item cost functions $\gamma_1, \ldots, \gamma_m \in 
		\Gamma$ of respective arities $r_1,\ldots,r_m$,
		\item $s_{j,1},\ldots, s_{j,r_j} \in \{1,\ldots,r\}$ for every $j \in \{1,\ldots,m\}$,
		\item positive rational numbers $\lambda_1,\ldots, \lambda_m \in \Q_{\geq0}$, and  $c\in \Q$,\end{itemize}such that for every $(x_1,\ldots,x_r)\in D^r$, it holds \[\rho(x_1,\ldots,x_r)=\sum_{j=1}^m\lambda_j\gamma_j(x_{s_{j,1}}, \ldots, x_{s_{j,r_j}})+ c.\]  We call \emph{affine convex cone}, or simply \emph{cone}, of $\g$,  the finite-valued language $\cone(\g)$ (with domain $D$) containing all cost functions from $\mathrm{wRel}_{D}$ that are affine conical combinations over $\g$.\end{definition}

In the next lemma, we show that solving $\vcsp(\cone(\g))$ is not harder than solving $\vcsp(\g)$.

\begin{lemma}\label{prop:strongexppwr}
	Let $\g$ be a finite-valued language with an arbitrary domain $D$  and with finite size, and let $\Delta$ be a finite-valued finite reduct  of $\cone( \g)$. Then there exists a polynomial-time many-one reduction from $\vcsp(\Delta)$ 
	to $\vcsp(\g)$.
	
\end{lemma}

\begin{proof}We claim that, for every set of variables $V:=\{v_1,\ldots,v_n\}$ and every finite sum  $\phi$ of cost functions from  $\Delta$ with at most $n$ free variables, there exists a  finite sum  $\phi'$ of cost 
	functions from $\g$  such that for every $u \in \Q$  there exists $u'\in \Q$ such that the following holds:\begin{gather}\nonumber\text{there exists }h \colon V \to D \text{ with cost } \phi(h(v_1),\ldots,h(v_n)\leq u \\\label{eq:polytimeeq1}\text{if, and only if,}\\   \nonumber \text{there exists } h' \colon V \to D\text{ with cost } \phi'(h'(v_1),\ldots, h'(v_n))\leq u'.\end{gather}%Let $\sigma$ be the (finite) signature of $\g$ and let $\tau \supseteq \sigma$ be the (finite) signature of $\Delta$. 
Fixed an arbitrary set of variables $V:=\{v_1,\ldots,v_n\}$ and an arbitrary  $u \in \mathbb Q$,  let us consider the instance $I=(V,\phi,u)$ of $\vcsp(\Delta)$ such that  \[\phi(v_1,\ldots,v_n):=\sum_{j=1}^m\hat{\gamma}_j(v_1^j,\ldots,v_{r_j}^j)\text{,}\]where $\hat{\gamma}_j \in \Delta$ is of arity $r_j$ and $v^j_i\in V$, for  $1\leq j \leq m$ and for $1\leq i \leq r_j$. 	By the definition of affine conical combination, for every $j \in \{1,\ldots,m\}$ there exist \begin{itemize}\item cost functions $\gamma_{j,1},\ldots,\gamma_{j,m_j} \in \g$ of respective arities $r_{j,1}, \ldots, r_{j,m_j}$,\item $s^{j}_{i,1},\ldots, s^{j}_{i,r_{j,i}} \in \{1,\ldots,r_j\}$ for every $i\in \{1,\ldots,m_j\}$,
		 \item numerical coefficients $\lambda_{j,1}, \ldots, \lambda_{j,m_j} \in \Q_{\geq 0}$, and $c_j \in \Q$ \end{itemize}such that for every $h\colon V \to D$ \[{\hat{\gamma}_j}(h(v^j_1),\ldots,h(v^j_{r_j}))=\sum_{i=1}^{m_j}\lambda_{j,i}\gamma_{j,i}(h(v_{s^j_{i,1}}^{j}), \ldots, h(v_{s^j_{i,r_{j,i}}}^{j}))+c_j.\]  We define an instance $I':=(V,\phi',u') $ of $\vcsp(\g)$ such that Condition~(\ref{eq:polytimeeq1}) holds,  $I'$ is computable in polynomial time in the size of $I=(V,\phi,u)$, and $\phi'$ does not depend on $u$. For every $1\leq j \leq m$  and for every $1\leq i\leq m_j$ there exist positive integers $\alpha_{j,i}$ and $\beta_{j,i}$ such that \[\lambda_{j,i}=\frac{\alpha_{j,i}}{\beta_{j,i}} \text{ and } \gcd(\alpha_{j,i}, \beta_{j,i})=1.\] For $1\leq j \leq m$ and for $1\leq i\leq m_{j}$, let us define \begin{align*}l_j&:=\lcm(\beta_{j,1},\ldots,\beta_{j,{m_j}}),\\l&:=\lcm(l_1,\ldots,l_m) \text{, and} \\ \mu_{j,i}&:=l \lambda_{j,i}.\end{align*}  We define \[\phi'(v_1,\ldots,v_n):=\sum_{j=1}^{m}\sum_{i=1}^{m_j}\mu_{j,i}\gamma_{j,i}(v_{s^j_{i,1}}^{j}, \ldots, v_{s^j_{i,r_{j,i}}}^{j}).\]Let $u':=l\left(u -\sum_{j=1}^mc_j\right)$. Therefore, for every $h \colon V \to D$ it holds that \[\phi(h(v_1),\ldots,h(v_n)< u\]if, and only if, \begin{align*}&\sum_{j=1}^{m} \sum_{i=1}^{m_j}\mu_{j,i}\gamma_{j,i}(h(v_{s^j_{i,1}}^{j}), \ldots, h(v_{s^j_{i,r_{j,i}}}^{j}))={\phi'}(h(v_1),\ldots,h(v_n))< u'.\qedhere \end{align*}\end{proof}

\begin{definition}\label{def:exprpwr2}
	Let $\g$ be a finite-valued language with an arbitrary domain $D$  and let $\mu 
	\colon D^r\to \QQ$ be a cost function. We say that $\mu$ is \emph{expressible in} $\g$ if there exist \begin{itemize}
		\item a set of variables $V:=\{v_1,\ldots,v_n\}$, 
		\item a tuple  $L:=(w_1,\ldots,w_r)$ of variables from $V$, and
		\item an $n$-ary cost function $\rho$ that is an affine conical combination over $\g$,
	\end{itemize}
	such that for every $(x_1,\ldots,x_r)\in D^r$, \[\mu(x_1,\ldots,x_r)=\inf_{\substack{h\colon V \to D\colon\\ h(w_i)=x_i, 1\leq i \leq r} }\rho(h(v_1),\ldots,h(v_n)).\]
	The \emph{expressive power} of $\g$ is the finite-valued language $\langle \g \rangle$ containing all cost functions with domain $D$ that are expressible in $\g$.
\end{definition}

In the next lemma, we show that solving $\vcsp(\langle \g \rangle)$ is not harder than solving $\vcsp(\cone( \g))$.

\begin{lemma}\label{lemma:weakexppr}
	Let $\g$ be a finite-valued language with finite size, and let $\Delta$ be a valued finite reduct of $\langle \g \rangle$. Then there exists a polynomial-time many-one reduction from $\vcsp(\Delta)$ to $\vcsp(\Theta)$, for some valued finite  reduct $\Theta$ of $\cone( \g)$.
\end{lemma}

\begin{proof}
	%	Let $\tau$ be the signature of $\Delta$ and let $\sigma\subseteq \tau$ be the  signature of $\cone(\g)$.
We claim that for every set of variables $V:=\{x_1,\ldots,x_r\}$ and every finite sum  
	$\phi$ of cost functions over $\Delta$ there exist a polynomial-time computable set of variables $V':=\{v_1,\ldots,v_n\}$, a tuple $L:=( w_1,\ldots,w_r)$ of variables from $V'$, and a finite sum  $\phi'$ of cost functions over $\cone(\g)$  such that for every $u \in \Q$  there exists $u'\in \Q$ such that  \begin{gather} \nonumber
	\text{there exists }h \colon V \to D \text{ with cost } \phi(h(x_1),\ldots,h(x_r))<u \\
	\label{eq:polytimeeq}\text{if, and only if, there exists } h' \colon V' \to D \\  \nonumber\text{ such that  }  h'(w_i)=h(x_i), 1\leq i \leq r \text{ with cost } \phi'(h'(v_1),\ldots, h(v_n))<u'.\end{gather}
	
	Let $I=(V,\phi,u)$ be an instance of $\vcsp(\Delta)$ such that $V:=\{x_1,\ldots,x_r\}$  and \[\phi(x_1,\ldots,x_r):=\sum_{j=1}^m\hat{\gamma}_j(v_1^j,\ldots,v_{r_j}^j)\]
	where $\hat{\gamma}_j \in \Delta$ and $v^j_i\in V$, for  $1\leq j \leq 
	m$ and for $1\leq i \leq r_j$. 	
	By the definition of expressive power~(Definition~\ref{def:exprpwr2}), for every $j \in \{1,\ldots,m\}$ there exist \begin{itemize}
		\item a  set of  new variables $W_j:=\{v^j_{r_j+1},\ldots,v^j_{n_j}\}$ such that $W_j \cap V=\emptyset$ and $W_j\cap W_{j'}=\emptyset$ for 
		every $j'\neq j$, $1\leq j'\leq m$ (it may be that $W_j=\emptyset$),
		\item an $n_j$-ary  cost function $\rho_j\in {\cone(\g)}$,
		\item $v^{j}_l \in V \cup W_j$, for $1\leq l\leq n_{j}$
	\end{itemize}
	such that for every $h\colon V \to D$ \[{\hat{\gamma}_j}(h(x^j_1),\ldots,h(x^j_{r_j}))=\inf_{\substack{h'\colon V\cup W_j \to D \colon \\ h'_{|V}=h}}\rho_{j}(h'(v_1^{j}), \ldots, h'(v_{{n_j}}^{j})).\]
	We define an instance $I':=(V',\phi',u')$ of $\vcsp(\cone(\g))$ such that Condition~(\ref{eq:polytimeeq}) holds,  $I'$ is computable in polynomial time in the size of $I=(V,\phi,u)$, and  $V'$ and $\phi'$ do not depend on $u$.
Set $V':=V \cup \bigcup_{j=1}^m W_j=\{v_1,\ldots,v_n\}$. We 
	define \[\phi'(v_1,\ldots,v_n):=\sum_{j=1}^{m}\rho_{j}(v^{j}_1,\ldots,v^{j}_{n_j}),\]
	and let $L:=\{x_1,\ldots,x_r\}$ be an enumeration of the variables in $V$.
	Let $u':=u$. Assume that $h \colon V \to D$ is such that \[\phi(h(x_1),\ldots,h(x_r))<u.\]
	Observe that, since  the sets of variables $W_j$ are mutually disjoint,
	\begin{align*}
	\phi(h(x_1), \ldots, h(x_r))&=&\sum_{j=1}^{m} \inf_{\substack{h'\colon V\cup W_j \to D \colon \\ h'_{|V}=h}}\rho_{j}(h'(v^{j}_1),\ldots,h'(v^{j}_{n_j}))\\
	&=&\inf_{\substack{h'\colon V' \to D \colon \\h'(w_i)=h(x_i), 1\leq 
			i \leq r}}{\phi'}(h'(v_1),\ldots,h'(v_n)).
	\end{align*}
	Hence, there exists %there exists an assignment $h \colon V \to \Q$ such that \[\phi^{\Delta}(h(v_1),\ldots,h(v_r))<u\] if, and only if, there exists 
	an assignment $h' \colon V' \to \Q$ such that  \[\phi'(h'(v_1),\ldots,h'(v_r),h'(v_{r+1}),\ldots, h'(v_n))<u.\]
	Let us assume now that there exists an assignment $h^\star\colon V' \to D$ such that \[\phi'(h^\star(v_1),\ldots,h^\star(v_n))<u,\] then it follows immediately  that \[\inf_{\substack{h'\colon V' \to D\colon \\ h'(w_i)=h^\star(x_i), 1\leq i \leq r}}\phi'(h'(v_1),\ldots,h'(v_n))< u, \]
	and, consequently, $\phi(h^\star(x_1),\ldots,h^\star(x_r))<u$. \qedhere
	
%	Vice versa, let us assume that there exists an assignment $h \colon V\to 	D$ such that $\phi(h(x_1),\ldots,h(x_n))<u$, % then \[\inf_{h'\colon V' \to D: h'(v_i)=h(v_i), 1\leq i \leq r}\phi'^{\Delta}(h'(v_1),\ldots,h'(v_n))< u'. \]	and let \[\varepsilon:=u' -\inf_{\substack{h'\colon V' \to D\colon \\ h'(w_i)=h(x_i), 1\leq i \leq r}}\phi'(h'(v_1),\ldots,h'(v_n)),\] then, by definition of the infimum, there exists $h^\star\colon V'\to D$ such that  $h^\star(w_i)=h(x_i)$ for $1\leq i \leq r$, and  \begin{align*}&\phi'(h^\star(v_1),\ldots,h^\star(v_n))<\inf_{\substack{h'\colon V' \to D\colon\\ h'(v_i)=h(x_i), 1\leq i \leq r}}\phi'	(h'(v_1),\ldots,h'(v_n))+\varepsilon=u'.\qedhere\end{align*} 
\end{proof}

\begin{proposition}\label{prop:weakexppr}
	Let $\g$ be a finite-valued language with an arbitrary domain $D$ having finite 
	size. Then for every  valued finite reduct $\Delta$ of $\langle \g \rangle$
	there exists a polynomial-time reduction from $\vcsp(\Delta)$ to $\vcsp(\g)$.
\end{proposition}

\begin{proof}
	To prove Proposition~\ref{prop:weakexppr}, it is enough to compose the polynomial-time reductions provided in Lemma~\ref{prop:strongexppwr} and in Lemma~\ref{lemma:weakexppr}. 
\end{proof}

\section{Local Expressive Power}\label{sect:exppwr2}

Informally, given a  finite-valued language $\g$ with domain $D$, the {\it local 
	expressive power} of $\g$  consists of  all cost functions whose restrictions to any finite subset $S$ of $D$ can be simulated by using the restrictions of cost functions from $\g$ to $S$.

\begin{definition}\label{definition:local.expressive.power.new} Let $D$ be a set and let us consider  $\Gamma \subseteq \mathrm{wRel}_{D}$. Let us 
	set \[\mathrm{M}_{k} (\Gamma) \defeq \left\{ (S,\gamma) \left| \, \gamma \in \Gamma, \, S \in D^{ \ar (\gamma) \times k}\right\} \right. \] for $k 
	\in \mathbb{N}$. We say that $\rho \colon D^{r} \to \mathbb{R}$  belongs to the \emph{local expressive power} of $\Gamma$ and write $\rho \in \mathrm{\ell Expr}(\Gamma)$ if, for every $\varepsilon > 0$ and every $k \in \mathbb{N}$, every $x^{1},\ldots,x^{k} \in D^{r}$, and every finite subset $\mathcal F \subseteq \mathcal{O}_{D}^{(k)}$ there exist $\lambda \in \mathbb{R}_{\geq 0}[\mathrm{M}_{k}(\Gamma)]$ and $c \in \mathbb{R}$ such that, for each $i \in \{ 1,\ldots,k \}$, \begin{equation}\label{eq:locexpr1}
	\left\lvert \rho(x^{i}) - \left( \sum_{(S,\gamma) \in \supp (\lambda)} \lambda(S,\gamma) \gamma \! \left(e_{i}^{(k)}(S)\right) + c \right) \right\rvert \, \leq \, \varepsilon 
	\end{equation} and, for each $f \in \mathcal F$, \begin{equation}\label{eq:locexpr2}
	\rho(f(x^{1},\ldots,x^{k})) \, \leq \, \sum_{(S,\gamma) \in \supp (\lambda)} \lambda(S,\gamma) \gamma (f(S)) + c + \varepsilon,
	\end{equation} where $e^{(k)}_i$ and $f$  are applied to the matrices $S$ componentwise, i.e.,  to their rows.\end{definition}

We equip $\mathrm{wRel}_{D}$ with the topology of pointwise convergence induced by the Euclidean topology on $\mathbb R$.

\begin{remark}\label{remark:closure} Let $D$ be a set and let $\g\subseteq \mathrm{wRel}_{D}$. \begin{enumerate}\item The  local expressive power $\mathrm{\ell Expr}(\Gamma)$ is a topologically closed subset of 
$\mathrm{wRel}_{D}$. To prove this, let $r \in \mathbb N$ and let $ \rho \in \mathrm{wRel}_{D}^{(r)}$ be in the closure of $\mathrm{\ell Expr}(\Gamma)$. Now, consider any $\varepsilon>0$, $r, k\in \mathbb N$, $x^1,\ldots,x^k \in D^r$, and a finite set $\mathcal F \subseteq \mathcal{O}_{D}^{(k)}$. As now $E \defeq \{x^1,\ldots,x^k\}\cup \{f(x^1,\ldots,x^k)\mid f \in \mathcal F\}$ is a finite subset of $D^{r}$, thus there exists $\tilde{\rho} \in \mathrm{\ell Expr}(\Gamma)$ such that   \[ \forall x \in E \colon \qquad \lvert\rho(x)-\tilde{\rho}(x)\rvert \, \leq \, \frac{\varepsilon }{2}.\] Since $\tilde{\rho} \in \mathrm{\ell Expr}(\Gamma)$, moreover there exist  $\lambda \in \mathbb{R}_{\geq 0}[\mathrm{M}_{k}(\Gamma)]$ and $c \in \mathbb{R}$ such that, for each $i \in \{ 1,\ldots,k \}$, \begin{align*}
&\left\lvert \rho(x^{i}) - \left( \sum\nolimits_{(S,\gamma) \in \supp (\lambda)} \lambda(S,\gamma) \gamma \! \left(e_{i}^{(k)}(S)\right) + c \right) \right\rvert \\
&\ \leq \, \left\lvert \rho(x^{i}) - \tilde{\rho}(x^{i}) \right\rvert + \left\lvert  \tilde{\rho}(x^{i}) - \left( \sum\nolimits_{(S,\gamma) \in \supp (\lambda)} \lambda(S,\gamma) \gamma \! \left(e_{i}^{(k)}(S)\right) + c \right) \right\rvert \\
&\ \leq \, \frac{\varepsilon }{2}+ \frac{\varepsilon }{2} \, = \, \varepsilon; 
\end{align*} and, for each $f \in \mathcal F$, \begin{align*}
 &\rho(f(x^{1},\ldots,x^{k}))-\left( \sum\nolimits_{(S,\gamma) \in \supp (\lambda)} \lambda(S,\gamma) \gamma (f(S)) + c \right)\\
 &\quad \leq \,  \left(\rho(f(x^{1},\ldots,x^{k}))- \tilde{\rho}(f(x^{1},\ldots,x^{k})) \right)\\
 &\qquad \quad +  \left(\tilde{\rho}(f(x^{1},\ldots,x^{k}))-\left(\sum\nolimits_{(S,\gamma) \in \supp (\lambda)} \lambda(S,\gamma) \gamma (f(S)) + c\right) \right)\\
 &\quad \leq \frac{\varepsilon }{2}+ \frac{\varepsilon }{2} \, = \,\varepsilon.
\end{align*}
This shows that $\rho \in  \mathrm{\ell Expr}(\Gamma)$.
\item A similar application of the triangle inequality, combined with the density of $\Q$ in $\R$, shows that in Definition~\ref{definition:local.expressive.power.new} the elements $\lambda \in \mathbb{R}_{\geq 0}[\mathrm{M}_{k}(\Gamma)]$ may  additionally be required to take only rational values.
\end{enumerate}
\end{remark}

\begin{proposition}\label{rem:exppwreqdef}
	Let $D$ be a set and let $\Gamma \subseteq \mathrm{wRel}_{D}$. Then $\overline{\langle \g \rangle}\subseteq \mathrm{\ell Expr}(\g)$, where $\overline{\langle \g \rangle}$ denotes the topological closure of $\langle \g \rangle$.
\end{proposition}

\begin{proof} By the first part of Remark~\ref{remark:closure} it suffices to check that ${\langle \g \rangle}\subseteq \mathrm{\ell Expr}(\g)$. Let us consider a cost function $\rho \colon D^r \to \Q$ such that $\rho \in \langle \g \rangle$. By the definition of expressive power (Definition~\ref{def:exprpwr2}) there exist 
	 a set of variables $V:=\{v_1,\ldots,v_n\}$, 
	a tuple $L:=( w_1,\ldots,w_r)$ of variables from $V$, 
cost functions  $\gamma_1, \ldots, \gamma_m $ from $\g$ with respective arities $r_1,\ldots,r_m$,  a column vector $(v^j_1,\ldots,v^j_{r_j})^{T}\in V^{r_j}$ for every $1\leq j\leq m$, $v^j$,
 positive rational numbers $\lambda_1,\ldots, \lambda_m$, and a rational number $c$,
	such that for every $(x_1,\ldots,x_r)\in D^r$ it holds \[\rho(x_1,\ldots,x_r)=\inf_{\substack{h\colon V \to D\colon\\ h(w_l)=x_l, 1\leq l \leq r} }\sum_{j=1}^m\lambda_j \gamma_j(h(v^j))+ c,\]
where the maps $h$ are applied componentwise.

Fix $\varepsilon>0$, a positive integer $k$, and   $x^1,\ldots,x^k\in D^r$.
	From the definition of $\rho$, it holds  that for every $i \in \{1,\ldots, k\}$, there exists   $s^i \colon V \to D$ with $s^i(w_l)=x^i_l$ for $1\leq l \leq r$ such that \begin{equation}\label{item:secondo} \left \lvert \rho(x^i)- 
	\left(\sum_{j=1}^m\lambda_j\gamma_j (s^i(v^j))+ c\right) \right\rvert \leq \varepsilon.\end{equation}
	For every $j\in \{1,\ldots,m\}$, let $S_j$ be the matrix from $D^{\ar(\gamma_j)\times k}$ with columns $s^1(v^j), \ldots, s^k(v^j)$ and  let $S_\rho \in D^{r\times k}$ be the matrix  whose columns are $x^1,\ldots x^k$.
	Let $\lambda \in \mathbb{R}_{\geq 0}[\mathrm{M}_{k}(\Gamma)]$ be such that \[\lambda(S,\gamma):=\begin{cases}
	\lambda_j & \text{if } (S,\gamma)=(S_j,\gamma_j)\\
	0 & \text{otherwise.}
	\end{cases}\]  
	With this notation,  we can immediately rewrite Condition~(\ref{item:secondo}) as \begin{itemize}
		\item[(i)]\label{item:secondolat}for  every $e^{(k)}_i \in \mathcal{J}_D^{(k)}$,   \[\left \lvert 
		\rho(e^{(k)}_i(x^1,\ldots,x^k))- \left(\sum_{(S,\gamma)\in \supp(\lambda)}\lambda(S,\gamma)\gamma(e^{(k)}_i(S))+ c\right) \right\rvert \leq \varepsilon.\] \end{itemize}
\noindent	Furthermore, the definition of $\rho$  implies, on the other hand, that for every tuple $y=(y_1,\ldots,y_r) \in D^r$ and  every $h \colon V \to D$ with $h(w_l)=y_l$ for $1\leq l 
		\leq r$ it holds \begin{equation} \label{item:primo}\rho(y)\leq \sum_{j=1}^m\lambda_j\gamma_j (h(v^j))+ c. \end{equation}
 For every $ f \in\mathcal O^{(k)}_D$, let $s \colon V \to D$ be the assignment defined by $ s(v):=f(s^1(v),\ldots, s^k(v)) \text{, for all } v \in V.$ Observe that \[ \begin{pmatrix} s(w_1)\\\vdots\\ s(w_r)\end{pmatrix}=f\left(\begin{pmatrix}
	s^1(w_1) & \ldots & s^k(w_1) \\
	\vdots &  & \vdots\\
	s^1(w_r) & \ldots & s^k(w_r) 
	\end{pmatrix}\right)=f(x^1,\ldots, x^k)\in D^r,\]
	and for every $j \in \{1,\ldots,m\}$ \[s(v^j):=\begin{pmatrix}s(v^j_1)\\\vdots\\s(v^j_{{\ar}(\gamma_j)})\end{pmatrix}=f\left(\begin{pmatrix}
	s^1(v^j_1) & \ldots & s^k(v^j_1) \\
	\vdots &  & \vdots\\
	s^1(v^j_{{\ar}(\gamma_j)}) & \ldots & s^k(v^j_{{\ar}(\gamma_j)}) 
	\end{pmatrix}\right)=f(S_j).\]
 By writing Condition~(\ref{item:primo}) for $y=f(x^1,\ldots,x^k)\in D^r$, we obtain the inequality $\rho(f(S_\rho))\leq \sum_{j=1}^{m}\lambda_j 
	\gamma_j(f(S_j))+c$, that is,
	\begin{itemize}
		\item[(ii)]\label{item:primolat}  for  every $f \in {\mathcal{O}}_D^{(k)}$,   \[ \rho(f(x^1,\ldots,x^k) )\leq \sum_{(S,\gamma)\in \supp(\lambda)}\lambda(S,\gamma)\gamma(f(S))+ c. \]
		\end{itemize}
	Observe that here we have proved a stronger condition than that required in Definition~\ref{definition:local.expressive.power.new}. Indeed, we have shown that the existence of the $\lambda$ and $c$ is universal (i.e.\ they are the same) for any $f \in {\mathcal{O}}_D^{(k)}$, and hence for every element of every arbitrary finite  $\mathcal F \subseteq  {\mathcal{O}}_D^{(k)}$.
Since Conditions~(i) and (ii) hold, $\rho \in \mathrm{\ell Expr}(\g)$ and $\langle\g\rangle\subseteq \mathrm{\ell Expr}(\g)$.  This completes the proof.\end{proof}

We now show  that for finite-valued languages with  finite domains, the local expressive power coincides with the topological closure of the expressive power. 

\begin{proposition}\label{prop:lexprfindom}
		Let $D$ be a finite set and consider  $\Gamma \subseteq \mathrm{wRel}_{D}$. Then $\overline{\langle \g \rangle}= \mathrm{\ell Expr}(\g)$, where $\overline{\langle \g \rangle}$ denotes the topological closure of $\langle \g \rangle$.
\end{proposition}

\begin{proof}
	We know that  $ \overline{\langle \g \rangle}\subseteq \mathrm{\ell Expr}(\g)$ (Proposition~\ref{rem:exppwreqdef}).
	Let us assume that $\rho \in \mathrm{\ell Expr}(\g)$. We claim that 
	for every   $\varepsilon >0$ there exists $\tilde{\rho}\in \langle \g \rangle$ such that \begin{displaymath}
	\forall x \in D^{r} \colon \quad \lvert\rho(x)-\tilde{\rho}(x)\rvert<2\varepsilon,
	\end{displaymath}  that is, $\rho \in \overline{\langle \g \rangle}$. 
	
	Let $k:=\lvert D \rvert^r$ and let $x^1,\ldots,x^k$ be an enumeration of all tuples from $D^r$. Because $D$ is a finite set,  the definition of local expressive power and the second part of Remark~\ref{remark:closure} asserts the existence of some rational-valued  $\lambda \in \R_{\geq0}\left [M_k(\g)\right]$ and $c \in \R$ such that  Conditions (\ref{eq:locexpr1}) and (\ref{eq:locexpr2}) of Definition~\ref{definition:local.expressive.power.new} hold for the finite set of operations $\mathcal F=\mathcal O^{(k)}_D$.
	Since $\lambda \in \R_{\geq0}\left [M_k(\g)\right]$ and  $\g$ has finite 
	domain and finite size, by definition there exist cost functions $\gamma_1, \ldots, \gamma_m \in \g$ (not necessarily pairwise distinct), and 
	matrices  $S_j \in D^{{\ar(\gamma_j)}\times k}$ for $1\leq j\leq m$ such that  $\supp(\lambda)=\{(\gamma_j,S_j)\mid 1\leq j \leq m\}$. Let us set $\lambda_j:=\lambda(\gamma_j, S_j)$ for $1\leq j \leq m$, and let $S_\rho \in 
	D^{r\times k}$ be the matrix with columns  $x^1, \ldots, x^k$. Then the following two conditions hold
	\begin{itemize}
		\item for all $e^{(k)}_i \in \mathcal J_D^{(k)}$,
		\begin{equation}	\label{cond:lexppwr1}  \left\lvert \rho(e^{(k)}_i(S_\rho))-\left(\sum_{j=1}^{m}\lambda_j \gamma_j(e^{(k)}_i(S_j))+c\right)\right\rvert \leq\varepsilon,\end{equation}
		\item for all $f \in \mathcal O_D^{(k)}$,
		\begin{equation}	\label{cond:lexppwr2}   \rho(f(S_\rho))-\left( \sum_{j=1}^{m}\lambda_j \gamma_j(f(S_j))+c\right)\leq \varepsilon.\end{equation}
	\end{itemize}
	Let us define the cost function $\tilde{\rho}\colon D^r \to \R$ such that 
	for every $i \in \{1,\ldots,k\}$, let \[\tilde{\rho}(x^i):=\min_{\substack{f \colon D^k \to D\colon\\ f(S_{\rho})=x^i}}\sum_{j=1}^{m}\lambda_j \gamma_j(f(S_j))+c,\] where the operations $f$ are applied componentwise, i.e., to the rows of the matrices $S_j$. The cost function $\tilde{\rho}$ is well defined 
because for every tuple $x=(x_1,\ldots,x_r) \in D^r$ there exists an index $i \in \{1,\ldots, k\}$ such that $(x_1,\ldots,x_r)=(x_1^i,\ldots,x_r^i)^{T}=x^i$.
	Observe that $\tilde{\rho}\in \langle \g \rangle$. In  order to see this, let us associate the rows of $S_\rho$ with fresh variables $\{v_1,\ldots,v_r\}$, and let us associate every row of every $S_j$ with a variable $v^j_l$, for  $1\leq l \leq \ar(\gamma_j)$. Let $V$ be the set of variables defined as \begin{align*}
	V & \, = \, \{v_1,\ldots,v_r,v_{r+1},\ldots,v_n\} \\
	&:= \, \{v_1,\ldots,v_r\}\cup \{v^j_l \mid 1\leq j \leq m \text{, } 1\leq l \leq \ar(\gamma_j)\},
	\end{align*}
	where two variables of $V$ are the same variable whenever they are associated with two rows that are equal as tuples of $D^k$.
Let us enumerate all the rows $a_1,\ldots, a_n$ of the matrices $S_{\rho}, S_1, \ldots, S_m$ and define the set  of their transposes $\mathcal A:=\{a_i^T\colon 1\leq i \leq n\}$. Then every assignment $s \colon V \to D$
	corresponds to a function $f \colon \mathcal A \to D$, and hence it holds that for every $x=(x_1\ldots,x_r) \in D^r$ \[\tilde{\rho}(x)=\min_{\substack{s \colon V \to D \colon \\ s(v_i)=x_i, \; 1\leq i \leq r}}\sum_{j=1}^{m}\lambda_j \gamma_j(s(v^j_1),\ldots,s(v^j_{\ar(\gamma_j)}))+c,\] that is, $\tilde{\rho} \in \langle \g \rangle$.
	
	Let $i \in \{1,\ldots,k\}$. We observe that, since $\mathcal F=\mathcal O^{(k)}_D$ there exists  $f^i \in \mathcal F$ such that $f^i(S_\rho)=x^i$ and $\tilde{\rho}(x^i)=\sum_{j=1}^{m}\lambda_j \gamma_j(f^i(S_j))+c$ (in other words, $f^i$ is the argument of the minimum defining $\tilde{\rho}(x^i)$); therefore,  by Condition~(\ref{cond:lexppwr2}) we obtain
	\begin{equation}\label{eq:finexppwrfin2}\rho(x^i)\leq \left(\sum_{j=1}^{m}\lambda_j \gamma_j(f^{i}(S_j))+c\right)+\varepsilon= \tilde{\rho}(x^i)+\varepsilon.\end{equation} On the other hand,  since $\mathcal{J}_{D}^{(k)} \subseteq \mathcal F$ and $e^{(k)}_i(S_\rho)=x^i$, from the definition of $\tilde{\rho}$ and from Condition~(\ref{cond:lexppwr1}) it follows that \begin{equation}\label{eq:finexppwrfin1}\tilde{\rho}(x^i)-\varepsilon \leq \left(\sum_{j=1}^{m}\lambda_j \gamma_j(e^{(k)}_i(S_j))+c\right)-\varepsilon\leq \rho(x^i).\end{equation}
	Therefore,  $\lvert\rho(x^i)-\tilde{\rho}(x^i)\rvert\leq\varepsilon<2\varepsilon$ for every 
	$i \in \{1,\ldots,k\}$, i.e., $\lvert\rho(x)-\tilde{\rho}(x)\rvert\leq\varepsilon<2\varepsilon$ for every $x \in D^r$. This proves the claim.
\end{proof}

\section{Weighted Polymorphisms}

\begin{definition} \label{def:weighting}Let $D$ be an arbitrary set.  For each $k \in \mathbb{N}$, we define the \emph{set of $k$-ary weightings on} $\mathcal O_D^{(k)}$ as \begin{displaymath}
	\mathrm{w}\mathcal{O}_{D}^{(k)} \, \defeq \, \left\{ \omega \in \mathbb{R}\!\left[ \mathcal{O}_{D}^{(k)} \right] \Bigg.\Bigg\rvert \sum_{f \in \supp (\omega)} \omega (f) = 0, \, \forall f \in  \mathcal{O}_{D}^{(k)}\setminus \mathcal{J}_{D}^{(k)} \colon \, \omega (f) \geq 0 \right\}.
	\end{displaymath} Moreover, let $\mathrm{w}\mathcal{O}_{D} \defeq \bigcup_{k \in \mathbb{N}} \mathrm{w}\mathcal{O}_D^{(k)}$. Given $\omega \in \mathrm{w}\mathcal{O}_{D}^{(k)}$ and $\rho \in \mathrm{wRel}_{D}^{(r)}$ with $k,r \in \mathbb{N}$, we say that $\omega$ is a \emph{weighted polymorphism} of $\rho$ (resp., $\rho$ is \emph{weight-improved} by $\omega$) if \begin{displaymath}
	\sum_{f \in \supp (\omega)} \omega (f) \rho (f(x^{1},\ldots,x^{k})) \, \leq \, 0
	\end{displaymath} for all $x^{1},\ldots,x^{k} \in D^{r}$. For $\Gamma \subseteq \mathrm{wRel}_{D}$, we define \begin{displaymath}
	\mathrm{wPol}(\Gamma) \, \defeq \, \left \{ \omega \in \mathrm{w}\mathcal{O}_{D} \mid \forall \gamma \in \Gamma \colon \, \omega \text{ is a weighted polymorphism of } \gamma \right\}.
	\end{displaymath} For $\Omega \subseteq \mathrm{w}\mathcal{O}_{D}$, we let \[\mathrm{Imp}(\Omega) \, \defeq \, \{ \gamma \in \mathrm{wRel}_{D} \mid \forall \omega \in \Omega \colon \, \gamma \text{ is weight-improved by } \omega \}.\] \end{definition}

Given an infinite  set  $D$, we aim at a description of the set  $\mathrm{Imp}(\mathrm{wPol}(\g))$ for $\Gamma \subseteq \mathrm{wRel}_{D}$.

Our main result is the following characterisation of the set $\mathrm{Imp}(\mathrm{wPol}(\Gamma))$ for finite-valued languages $\Gamma$ with  infinite domains.

\begin{theorem}\label{theorem:main} Let $D$ be an arbitrary infinite set 
	and let $\Gamma \subseteq \mathrm{wRel}_{D}$. Then \[\mathrm{Imp}(\mathrm{wPol}(\Gamma)) = \mathrm{\ell Expr}(\Gamma).\] \end{theorem}

\begin{remark}
	We remark that Theorem~\ref{theorem:main} is consistent with the characterisation of $\mathrm{Imp}(\mathrm{wPol}(\Gamma))$ for finite-domain finite-valued languages~\cite{FullaZivny,CohenCooperCreedJeavonsZivny}. Indeed, for a finite-domain finite-valued language $\g$ it holds that \[\mathrm{Imp}(\mathrm{wPol}(\Gamma)) = \overline{\langle\Gamma\rangle}=\mathrm{\ell Expr}(\Gamma),\]
	where the first inequality follows from~\cite[Theorem~3.3]{FullaZivny}, and the second inequality  follows from Proposition~\ref{rem:exppwreqdef} applied to the specific case of a finite-domain finite-valued language.
\end{remark}

\begin{corollary}\label{lemma:expsubsetimpfpol}Let $\g$ be a finite-valued language with an arbitrary  domain $D$ and values in $\R$.
	Then every weighted polymorphism of $\g$ improves every cost function that is 
	expressible in $\g$, that is, \[\langle \g \rangle \subseteq \imp \left( \wpol (\g)\right). \]
\end{corollary}

\begin{proof}
This follows immediately from a combination of Theorem~\ref{theorem:main} and Proposition~\ref{rem:exppwreqdef}. Below we include an alternative direct proof, illustrating that the crucial inclusion established by our Theorem~\ref{theorem:main} is ($\subseteq$).

 Let $\omega$ be a $k$-ary weighted polymorphism of $\g$, and let $\rho$ be a cost function expressible in $\g$. Let $r \defeq \ar (\rho)$. Consider any $x^1,\ldots,x^k \in D^r$. We aim to prove that \begin{equation}\label{formula:goal}
\sum\nolimits_{g \in \supp(\omega)}\omega(g)\rho(g(x^1,\ldots, x^k)) \, \leq \, 0 ,
\end{equation} Since~\eqref{formula:goal} holds trivially if $\omega=0$, we may and will assume that $\omega \ne 0$. Let \begin{displaymath}
\supp^{+}(\omega) \, \defeq \, \{g \in \supp(\omega)\mid \omega(g)>0\}
\end{displaymath} and $\supp^{-}(\omega) \defeq \supp(\omega)\setminus\supp^{+}(\omega)$. Note that $\supp^{-}(\omega) \subseteq \mathcal{J}_{D}^{(k)}$ by the definition of a weighting (Definition~\ref{def:weighting}). As $\omega \ne \emptyset$, we have $\supp^{-}(\omega)\neq \emptyset$. In order to prove~\eqref{formula:goal}, it suffices to show that \begin{equation}\label{eq:goal2}
\forall \varepsilon \in \R_{>0} \colon \qquad \sum\nolimits_{g \in \supp(\omega)}\omega(g)\rho(g(x^1,\ldots, x^k)) \, \leq \, \varepsilon .
\end{equation} For this purpose, let $\varepsilon \in \R_{>0}$. We define \begin{displaymath}
\varepsilon^\star \, \defeq \, -\frac{\varepsilon}{\sum_{e \in \supp^{-}(\omega)}\omega(e)}
\end{displaymath} and observe that $\varepsilon^\star>0$. As $\rho$ is expressible in $\g$ (cf.~Definition \ref{def:exprpwr2}), there exist a finite set of variables $V:=\{v_1,\ldots, v_n\}$, a tuple $L\defeq ( w_1,\ldots,w_r)$ of variables from $V$,  a finite sequence of cost functions $\gamma_1,\ldots, \gamma_m \in \g$, tuples $v^1 \in V^{\ar(\gamma_1)},\ldots,v^m \in V^{\ar(\gamma_m)}$, coefficients $\lambda_1,\ldots, \lambda_m \in \R_{\geq 0}$, and a constant $c \in \R$ such that, for every $z \in D^r$, \begin{displaymath}
\rho(z) \, = \, \rho(z_1,\ldots,z_r)\, = \, \inf_{\substack{h\colon V \to D\colon \\ h(w_i)=z_i, \, 1\leq i \leq r }}\sum\nolimits_{j=1}^m\lambda_j \gamma_j (h(v^j))+c.
\end{displaymath} In particular, for every $g \in \mathcal{O}_{D}^{(k)}$, \begin{equation}\label{eq:defrho2}
\rho(g(x^{1},\ldots,x^{k})) \, = \, \inf_{\substack{h\colon V \to D\colon \\ h(w_i)=g(x_{i}^{1},\ldots,x_{i}^{k}), \, 1\leq i \leq r }}\sum\nolimits_{j=1}^m\lambda_j \gamma_j (h(v^j))+c.
\end{equation} Thus, by~\eqref{eq:defrho2}, for each $l \in \{ 1,\ldots,k\}$, there exists an assignment $h^{l} \colon V \to D$ such that $h^{l}(w_1)=x_{1}^{l}, \ldots, h^{l}(w_r)=x_{r}^{l}$ and \begin{align}\label{eq:expsubsetimpfpol2}
\sum\nolimits_{j=1}^m\lambda_j \gamma_j(h^{l}(v^j))+c \, \leq \, \rho(x^{l}) + \varepsilon^\star.
\end{align} For each $g \in \mathcal{O}_{D}^{(k)}$, since the assignment $h^{\ast}_{g} \colon V \to D, \, u \mapsto g(h^1(u), \ldots,h^k(u))$ satisfies \begin{displaymath}
\forall i \in \{ 1,\ldots,r \} \colon \quad h^{\ast}_{g}(w_{i}) \, = \, g(h^1(w_{i}),\ldots,h^k(w_{i})) \, = \, g(x_{i}^{1},\ldots,x_{i}^{k}) ,
\end{displaymath} assertion~\eqref{eq:defrho2} moreover entails that \begin{equation}\label{eq:expsubsetimpfpol1}
\rho(g(x^{1},\ldots,x^{k})) \, \leq \, \sum\nolimits_{j=1}^m \lambda_j \gamma_j(g(h^1(v^j), \ldots,h^k(v^j)) )+c .
\end{equation} Since $\omega$ is a weighted polymorphism of $\g$, we know that $\omega$ weight-improves each of $\gamma_1,\ldots, \gamma_m$ (cf.~Definition~\ref{def:weighting}). In particular, \begin{equation}\label{eq:impro}
\forall j\in \{1,\ldots,m\} \colon \quad \sum\nolimits_{g \in \supp(\omega)}\omega(g)\gamma_j(g(h^1(v^j)), \ldots,h^k(v^j)) \, \leq \, 0 .
\end{equation} We now conclude that
\begin{align*}
&\sum\nolimits_{g \in \supp(\omega)}\omega(g)\rho(g(x^1,\ldots, x^k))\\
& \quad= \, \sum\nolimits_{g \in \supp(\omega)\setminus \mathcal J_D^{(k)}}\omega(g)\rho(g(x^1,\ldots, x^k)) +\sum\nolimits_{l\colon e_l^{(k)} \in \supp^{-}(\omega)}\omega\!\left(e_l^{(k)} \right)\!\rho(x^l) \\   
&\quad \stackrel{\eqref{eq:expsubsetimpfpol1}+\eqref{eq:expsubsetimpfpol2}}{\leq} \, \sum\nolimits_{g \in \supp^+(\omega) }\omega(g)\!\left(\sum\nolimits_{j=1}^m \lambda_j \gamma_j(g(h^1(v^j), \ldots,h^k(v^j)) )+c  \right)\\
& \qquad \qquad \qquad +\sum\nolimits_{l \colon e_l^{(k)} \in \supp^{-}(\omega)}\omega\!\left(e_l^{(k)} \right)\!\left(\sum\nolimits_{j=1}^m\lambda_j \gamma_j(h^l(v^j))+c-\varepsilon^\star  \right) \\
&   \quad = \,  \sum\nolimits_{j=1}^m \lambda_j \!\left(\sum\nolimits_{g \in \supp(\omega)}\omega(g)\gamma_j(g(h^1(v^j), \ldots,h^k(v^j)) \right)\\
& \qquad \qquad \qquad + \left(\sum\nolimits_{g \in \supp(\omega)}\omega(g)\right)\! c -\left(\sum\nolimits_{e_l^{(k)} \in \supp^{-}(\omega)}\omega\!\left(e_l^{(k)} \right)\!\right)\!\varepsilon^\star\\
& \quad \leq \,  \sum\nolimits_{j=1}^m \lambda_j \!\left( \sum\nolimits_{g \in \supp(\omega)}\omega(g)\gamma_j(g(h^1(v^j)), \ldots,h^k(v^j)) \right)+\varepsilon \, \stackrel{\eqref{eq:impro}}{\leq} \, \varepsilon,
\end{align*} where the second to last inequality follows from $\sum_{g \in \supp(\omega)}\omega(g)=0$ and the choice of $\varepsilon^\star$, and the last estimate is using non-negativity of $\lambda_1,\ldots, \lambda_m$, too. This proves~\eqref{eq:goal2} and hence completes the argument. \end{proof}

Before proceeding to the proof of Theorem~\ref{theorem:main}, we recollect some basic functional-analytic tools in the following section.

\section{Locally Convex Spaces}\label{section:constructions}
In this section, we compile some background material on locally convex
topological vector spaces relevant for the proof of
Theorem~\ref{theorem:main}.

Throughout this note, the term \emph{vector space} will always mean a
\emph{vector space over the field $\mathbb{R}$}. To clarify some basic
geometric terminology, let $X$ be a vector space. A subset $S \subseteq
X$ is said to be \begin{itemize}
	\item[---] a \emph{cone} if $\lambda x \in S$ for all $x \in S$ and
	$\lambda \in \R_{>0}$,
	\item[---] \emph{convex} if $\lambda x + \mu y \in S$ for all $x,y
	\in S$ and $\lambda , \mu \in \R_{\geq 0}$ with $\lambda + \mu = 1$,
	\item[---] \emph{absorbent} if \begin{displaymath}
	\forall x \in X \ \exists \alpha \in \R_{>0} \ \forall \lambda \in
	\R \colon \quad \vert \lambda \vert \geq \alpha \ \Longrightarrow \ x
	\in \lambda S ,
	\end{displaymath}
	\item[---] \emph{symmetric} if $-S = S$. \end{itemize}

A \emph{topological vector space} is a vector space $X$ equipped with a
topology such that both $X \times X \to X, \, (x,y) \mapsto x+y$ and $\R
\times X \to X, \, (\lambda,x) \mapsto \lambda x$ are continuous maps.
In order to agree on some further terminology and notation, let $X$ be a
topological vector space. As usual, we say that $X$ is \emph{locally
	convex} if $0 \in X$ admits a neighbourhood base consisting of convex
subsets of $X$. We let $X^{\ast}$ denote the \emph{topological dual} of
$X$, i.e., the topological vector space of all continuous linear forms
on $X$ endowed with the weak-$*$ topology, which is the initial topology
on $X^{\ast}$ generated by the maps of the form \begin{displaymath}
X^{\ast} \! \, \longrightarrow \, \mathbb{R}, \quad \phi \,
\longmapsto \, \phi(x) \qquad (x \in X) .
\end{displaymath} It is well known that $X^{\ast}$ is a locally convex
Hausdorff topological vector space (cf.~\cite[II.42--43,
\S6.2]{bourbaki}). For a subset $S \subseteq X$, we define
\begin{displaymath}
S^{+} \! \, \defeq \, \{ \phi \in X^{\ast} \mid \forall x \in S
\colon \, \phi(x) \geq 0 \} .
\end{displaymath} Furthermore, if $Y$ is another topological vector
space and $A \colon X \to Y$ is linear and continuous, then we consider
the continuous linear map \begin{displaymath}
A^{\ast} \colon \, Y^{\ast}\! \, \longrightarrow \, X^{\ast}, \quad
\psi \, \longmapsto \, \psi \circ A .
\end{displaymath} We now recall a fairly general version of Farkas' lemma.

\begin{lemma}[abstract core of Farkas' lemma; \cite{swartz},
	Lemma~2.1]\label{lemma:farkas} Let $X$ and $Y$ be locally convex
	Hausdorff topological vector spaces. If $S$ is a closed convex cone in
	$Y$ and $A \colon X \to Y$ is linear and continuous, then
	\[(A^{-1}(S))^{+} \, = \, \overline{A^{\ast}(S^{+})}.\] \end{lemma}

The following more familiar version of Farkas' lemma is a mere
reformulation of Lemma~\ref{lemma:farkas}.

\begin{theorem}[Farkas' lemma]\label{theorem:farkas} Let $X$ and $Y$ be
	locally convex Hausdorff topological vector spaces and let $A \colon X
	\to Y$ be linear and continuous. Let $S$ be a closed convex cone in $Y$
	and let $\phi \in X^{\ast}$. The following are equivalent. \begin{enumerate}
		\item[$(1)$] $\forall x \in X \colon \ A(x) \in S \,
		\Longrightarrow \,
		\phi(x) \geq 0$.
		\item[$(2)$] $\phi \in \overline{A^{\ast}(S^{+})}$.
\end{enumerate} \end{theorem}

\begin{proof} We deduce that \begin{align*}
	\phi \in \overline{A^{\ast}(S^{+})} \quad &
	\stackrel{\ref{lemma:farkas}}{\Longleftrightarrow} \quad \phi \in
	(A^{-1}(S))^{+} \\
	& \Longleftrightarrow \quad \forall x \in A^{-1}(S) \colon \ \phi
	(x) \geq 0 \\
	& \Longleftrightarrow \quad \forall x \in X \colon \ \left( A(x) \in
	S \, \Longrightarrow \, \phi (x) \geq 0 \right) . \qedhere
	\end{align*} \end{proof}

For our purposes, we will need a specific variation of
Theorem~\ref{theorem:farkas}, whose proof involves the following basic
observation.

\begin{remark}\label{remark:pointed.cone} Let $C$ be any non-empty
	closed cone in a topological vector space $X$. Then $0 \in C$. Indeed,
	picking any $x \in X$, we see that \begin{displaymath}
	0 \, = \, \lim\nolimits_{\R_{>0}\ni \lambda \to 0} \lambda x \,
	\in \, C .
	\end{displaymath} \end{remark}

\begin{corollary}\label{corollary:farkas} Let $X,Y,Z$ be locally convex
	Hausdorff topological vector spaces, and let $A \colon X \to Y$ and $B
	\colon X \to Z$ be linear and continuous. Moreover, let $S \subseteq Y$
	and $T \subseteq Z$ be non-empty closed convex cones and let $\phi \in
	X^{\ast}$. The following are equivalent. \begin{enumerate}
		\item[$(1)$] $\forall x \in X \colon \ (A(x) \in S \wedge B(x) \in
		T) \, \Longrightarrow \, \phi(x) \geq 0$.
		\item[$(2)$] $\phi \in \overline{\{ A^{\ast}(\mu) + B^{\ast}(\nu)
			\mid \mu \in S^{+}, \, \nu \in T^{+} \}}$.
\end{enumerate} \end{corollary}

\begin{proof} Let us consider the locally convex Hausdorff topological
	vector spaces $\tilde{X} \defeq X$ and $\tilde{Y} \defeq Y \times Z$,
	the continuous linear maps $\tilde{\phi} \defeq \phi$ and
	\begin{displaymath}
	\tilde{A} \colon \, \tilde{X} \, \longrightarrow \, \tilde{Y}, \quad
	x \, \longmapsto \, (A(x),B(x)) ,
	\end{displaymath} and the closed convex cone $\tilde{S} \defeq S \times
	T \subseteq \tilde{Y}$. It is well known that the mapping $\Psi \colon
	Y^{\ast} \times Z^{\ast} \to (Y \times Z)^{\ast}$ given by
	\begin{displaymath}
	\Psi (\mu, \nu) (y,z) \, \defeq \, \mu(y) + \nu (z) \qquad (\mu
	\in Y^{\ast}, \, \nu \in Z^{\ast}, \, y \in Y, \, z \in Z)
	\end{displaymath} constitutes an isomorphism of topological vector
	spaces. Note that, if $\mu \in Y^{\ast}$ and $\nu \in Z^{\ast}$, then
	\begin{align*}
	\tilde{A}^{\ast}(\Psi(\mu,\nu))(x) \, &= \,
	\Psi(\mu,\nu)(\tilde{A}(x)) \, = \, \Psi(\mu,\nu)(A(x),B(x)) \\
	& = \, \mu (A(x)) + \nu(B(x)) \, = \, A^{\ast}(\mu)(x) +
	B^{\ast}(\nu)(x) \\
	& = \, (A^{\ast}(\mu) + B^{\ast}(\nu))(x)
	\end{align*} for every $x \in X$, i.e., $\tilde{A}^{\ast}(\Psi
	(\mu,\nu)) = A^{\ast}(\mu) + B^{\ast}(\nu)$. Furthermore,
	Remark~\ref{remark:pointed.cone} asserts that $0_{Y} \in S$ and $0_{Z}
	\in T$. Form this and the definition of $\Psi$, we now conclude that
	\begin{displaymath}
	\forall \mu \in Y^{\ast} \ \forall \nu \in Z^{\ast} \colon \quad
	\, \Psi(\mu,\nu) \in (S \times T)^{+} \! \ \, \Longleftrightarrow \ \,
	\mu \in S^{+} \! \, \wedge \, \nu \in T^{+} .
	\end{displaymath} Consequently, \begin{align*}
	\tilde{A}^{\ast}(\tilde{S}^{+}) \, &= \, \left\{
	\tilde{A}^{\ast}(\Psi(\mu,\nu)) \left\vert \, \mu \in Y^{\ast}, \, \nu
	\in Z^{\ast}, \, \Psi(\mu,\nu) \in \tilde{S}^{+} \right\} \right. \\
	& = \, \left\{ A^{\ast}(\mu) + B^{\ast}(\nu) \left\vert \, \mu
	\in S^{+}, \, \nu \in T^{+} \right\} . \right.
	\end{align*} In the light of these observations, the assertion of
	Corollary~\ref{corollary:farkas} for the objects $X,Y,Z,\phi,A,B,S,T$ is
	an immediate consequence of Theorem~\ref{theorem:farkas} applied to
	$\tilde{X},\tilde{Y},\tilde{\phi},\tilde{A},\tilde{S}$. \end{proof}

In order to the accommodate the application of
Corollary~\ref{corollary:farkas} in the proof of
Theorem~\ref{theorem:main}, we now recall two standard constructions for
locally convex topological vector spaces.

\begin{proposition}[\cite{bourbaki}, II.25,
	\S4.2(2)]\label{proposition:bourbaki} Let $X$ be a vector space.
	\begin{itemize}
		\item[$(1)$] The set of all absorbent, symmetric, convex subsets of
		$X$ is a neighbourhood basis at $0$ for a unique vector space topology
		$\tau_{X}$ on $X$.
		\item[$(2)$] $\tau_{X}$ is the finest locally convex vector space
		topology on $X$.
		\item[$(3)$] $\tau_{X}$ is Hausdorff.
		\item[$(4)$] Any linear map from $X$ to another locally convex
		topological vector space is continuous with respect to $\tau_{X}$.
\end{itemize} \end{proposition}

\begin{proof} A proof is to be found in~\cite[II.25--26,
	\S4.2(2)]{bourbaki} (for an extension, see also~\cite[II.27, \S4.4,
	Proposition~5]{bourbaki}). \end{proof}

We proceed to another well-known construction of locally convex
topological vector spaces. If $X$ and $Y$ are two vector spaces and
$\beta \colon X \times Y \to \mathbb{R}$ is a bilinear map, then we let
$\sigma_{\beta}(X,Y)$ denote the initial topology on $X$ generated by
all maps of the form $X \to \mathbb{R}, \, x \mapsto \beta (x,y)$ with
$y \in Y$.

\begin{remark}[see~\cite{bourbaki}, II.42--43,
	\S6.2]\label{remark:weak.topologies} Let $X$ and $Y$ be two vector
	spaces and let $\beta \colon X \times Y \to \mathbb{R}$ be any bilinear
	map. The following hold. \begin{itemize}
		\item[$(1)$] $(X,\sigma_{\beta}(X,Y))$ is a locally convex
		topological vector space.
		\item[$(2)$] $\sigma_{\beta}(X,Y)$ is Hausdorff if and only if
		$\beta$ \emph{separates} $X$, i.e., \begin{displaymath}
		\qquad \forall x \in X \colon \quad (\forall y \in Y \colon \, \beta
		(x,y) = 0) \ \Longrightarrow \ x=0.
		\end{displaymath}
		\item[$(3)$] If $\beta$ \emph{separates} $Y$, i.e., \begin{displaymath}
		\qquad \forall y \in Y \colon \quad (\forall x \in X \colon \, \beta
		(x,y) = 0) \ \Longrightarrow \ y=0,
		\end{displaymath} then \begin{displaymath}
		Y \, \longrightarrow \, (X,\sigma_{\beta}(X,Y))^{\ast}, \quad y \,
		\longmapsto \, (x \mapsto \beta (x,y))
		\end{displaymath} is an isomorphism of vector spaces.
\end{itemize} \end{remark}

\section{Proof of Theorem~\ref{theorem:main}}

Everything is in place for proving
Theorem~\ref{theorem:main}. Our argument runs by an application of
Farkas' lemma in the form of Corollary~\ref{corollary:farkas}.

\begin{proof}[Proof of Theorem~\ref{theorem:main}]  ($\supseteq$): Let $r \in
	\mathbb{N}$ and $\rho \colon D^{r} \to \mathbb{R}$. For contradiction,
	let us assume that $\rho \in \mathrm{\ell Expr}(\Gamma) \setminus
	\mathrm{Imp}(\mathrm{wPol}(\Gamma))$. Then there exist $k \in
	\mathbb{N}$ and $\omega \in \mathrm{w}\mathcal{O}_{D}^{(k)} \cap
	\mathrm{wPol}(\Gamma)$ such that $\rho$ is not weight-improved by
	$\omega$, that is, we can find  $x^{1},\ldots,x^{k} \in D^{r}$ such that
	$\varepsilon \defeq \sum_{f \in \supp (\omega)} \omega (f) \rho
	(f(x^{1},\ldots,x^{k})) > 0.$ Let $\mathcal F \defeq \supp (\omega)$ and
	$L \defeq (\vert \mathcal F \vert + 1)(\sup_{f \in \mathcal F} \vert
	\omega(f) \vert + 1)$. Since $\rho \in \mathrm{\ell Expr}(\Gamma)$,
	there exist $\lambda \in \mathbb{R}_{\geq 0}[\mathrm{M}_{k}(\Gamma)]$
	and $c \in \mathbb{R}$ such that, for each $i \in \{ 1,\ldots,k \}$,
	\begin{displaymath}
	\left\lvert \rho(x^{i}) - \left( \sum\nolimits_{(S,\gamma) \in \supp
		(\lambda)} \lambda(S,\gamma) \gamma \! \left(e_{i}^{(k)}(S)\right) + c
	\right) \right\rvert \, \leq \, \frac{\varepsilon}{2L}
	\end{displaymath} and, for each $f \in \mathcal F$, \begin{displaymath}
	\rho(f(x^{1},\ldots,x^{k})) \, \leq \, \sum\nolimits_{(S,\gamma) \in
		\supp (\lambda)} \lambda(S,\gamma) \gamma (f(S)) + c +
	\frac{\varepsilon}{2L}.
	\end{displaymath} Now, if $f \in \mathcal F\setminus
	\mathcal{J}_{D}^{(k)}$, then $\omega(f) \geq 0$, and therefore
	\begin{align*}
	\omega (f) \rho &(f(x^{1},\ldots,x^{k})) \\
	&  \leq \, \omega (f)\! \left( \sum\nolimits_{(S,\gamma) \in \supp
		(\lambda)} \lambda(S,\gamma) \gamma (f(S)) + c + \frac{\varepsilon}{2L} 
	\right) \\
	& \leq  \, \omega (f)\! \left( \sum\nolimits_{(S,\gamma) \in \supp
		(\lambda)} \lambda(S,\gamma) \gamma (f(S)) + c \right) \! +
	\frac{\varepsilon}{2(\vert \mathcal{F} \vert + 1)}.
	\end{align*} On the other hand, if $f \in \mathcal F \cap
	\mathcal{J}_{D}^{(k)}$, then $f = e_{i}^{(k)}$ for some $i \in \{
	1,\ldots,k \}$, and
	hence \begin{align*}
	& \omega (f) \rho (f(x^{1},\ldots,x^{k})) \, = \, \omega (f) \rho
	(x^{i}) \\
	& \hspace{20mm} \leq \, \omega (f)\! \left(
	\sum\nolimits_{(S,\gamma) \in \supp (\lambda)} \lambda(S,\gamma) \gamma
	\! \left(e_{i}^{(k)}(S)\right) + c \right) + \frac{\varepsilon}{2(\vert
		\mathcal{F} \vert + 1)} \\
	& \hspace{20mm} = \, \omega (f) \!\left( \sum\nolimits_{(S,\gamma)
		\in \supp (\lambda)} \lambda(S,\gamma) \gamma (f(S)) + c \right) +
	\frac{\varepsilon}{2(\vert \mathcal{F} \vert + 1)}.
	\end{align*} We now conclude that \begin{align*}
	& \sum\nolimits_{f \in \supp (\omega)} \omega (f) \rho
	(f(x^{1},\ldots,x^{k})) \\
	& \ \ \leq \, \sum\nolimits_{f \in \supp (\omega)} \! \left( \omega
	(f)\! \left( \sum\nolimits_{(S,\gamma) \in \supp (\lambda)}
	\lambda(S,\gamma) \gamma (f(S)) + c \right) \! +
	\frac{\varepsilon}{2(\vert \mathcal{F} \vert +1)} \right) \\
	& \ \ \leq \, \frac{\varepsilon}{2} + \sum\nolimits_{f \in \supp
		(\omega)} \omega (f) \! \left( \sum\nolimits_{(S,\gamma) \in \supp
		(\lambda)} \lambda(S,\gamma) \gamma (f(S)) + c \right) \\
	& \ \ = \, \frac{\varepsilon}{2} + \sum\nolimits_{(S,\gamma) \in
		\supp (\lambda)} \lambda(S,\gamma) \! \left( \sum\nolimits_{f \in \supp
		(\omega)} \omega (f) \gamma (f(S)) \right)\! \, \leq \,
	\frac{\varepsilon}{2},
	\end{align*} thus arriving at the desired contradiction.
	
($\subseteq$): Let us consider some $r \in \mathbb{N}$ as well as a cost
	function $\rho \colon D^{r} \to \mathbb{R}$, and suppose that $\rho \in
	\mathrm{Imp}(\mathrm{wPol}(\Gamma))$. We will show that $\rho \in
	\mathrm{\ell Expr}(\Gamma)$. For this purpose, consider any $k \in
	\mathbb{N}$ and $x^{1},\ldots,x^{k} \in D^{r}$. Let $M \defeq
	\mathrm{M}_{k} (\Gamma)$ and put $N \defeq \mathcal{O}_{D}^{(k)}$.
	Preparing an application of Corollary~\ref{corollary:farkas}, let us
	consider the locally convex Hausdorff topological vector space $Z \defeq
	\left(\mathbb{R}[N],\tau_{\mathbb{R}[N]}\right)$, as defined in
	Proposition~\ref{proposition:bourbaki}. The linear subspace\footnote{It
		follows by Proposition~\ref{proposition:bourbaki}(4) that $X$ is closed
		in $Z$.} \begin{displaymath}
	X \, \defeq \, \left\{ \omega \in \mathbb{R}[N] \left\vert \,
	\sum\nolimits_{f \in \supp (\omega)} \omega(f) = 0 \right\} \right.
	\end{displaymath} of $Z$, endowed with the subspace topology
	inherited from $Z$, constitutes a locally convex Hausdorff topological
	vector space, too. Furthermore, we note that the bilinear map
	\begin{displaymath}
	\beta \colon \, \mathbb{R}^{M}\! \times \mathbb{R}[M] \,
	\longrightarrow \, \mathbb{R}, \quad (x,\lambda) \, \longmapsto \,
	\sum\nolimits_{(S,\gamma) \in \supp (\lambda)} \lambda(S,\gamma) x(S,\gamma)
	\end{displaymath} separates both $\mathbb{R}^{M}$ and
	$\mathbb{R}[M]$. Therefore, invoking
	Remark~\ref{remark:weak.topologies}(1)+(2), we see that $Y \defeq
	(\mathbb{R}^{M},\sigma_{\beta} (\mathbb{R}^{M},\mathbb{R}[M]))$ is
	locally convex Hausdorff topological\footnote{Of course, $\sigma_{\beta}
		(\mathbb{R}^{M},\mathbb{R}[M])$ coincides with the product topology on
		$\mathbb{R}^{M}$} vector space. Due to
	Proposition~\ref{proposition:bourbaki}(4), the three linear maps
	\begin{displaymath}
	\phi_{0} \colon \, Z \, \longrightarrow \, \mathbb{R}, \quad
	\omega \, \longmapsto \, \sum\nolimits_{f \in \supp (\omega)} \omega (f)
	\rho (f(x^{1},\ldots,x^{k})) ,
	\end{displaymath} $A_{0} \colon Z \to Y$ defined by \begin{displaymath}
	A_{0}(\omega)(S,\gamma) \, \defeq \, \sum\nolimits_{f \in \supp
		(\omega)} \omega (f) \gamma (f(S)) \qquad \left(\omega \in \R[N], \,
	(S,\gamma) \in M\right) ,
	\end{displaymath} and $B_{0} \colon Z \to Z$ defined by
	\begin{displaymath}
	B_{0}(\omega)(f) \, \defeq \, \begin{cases}
	\, \omega (f) & \text{if } f \notin \mathcal{J}_{D}^{(k)}, \\
	\, 0 & \text{if } f \in \mathcal{J}_{D}^{(k)}
	\end{cases} \qquad \left( \omega \in \R[N], \, f \in N\right)
	\end{displaymath} are continuous. Hence, \begin{displaymath}
	\phi \defeq \phi_{0}\vert_{X} \colon X \longrightarrow \R, \quad
	A \defeq A_{0}\vert_{X} \colon X \longrightarrow Y, \quad B \defeq
	B_{0}\vert_{X} \colon X \longrightarrow Z
	\end{displaymath} are continuous linear maps. We note that
	\begin{align*}
	S \, &\defeq \, \left. \! \left\{ x \in \mathbb{R}^{M}
	\right\vert \forall (S,\gamma) \in M \colon \, x(S,\gamma) \geq 0
	\right\} \\
	&= \, \left. \! \left\{ x \in \mathbb{R}^{M} \right\vert \forall
	(S,\gamma) \in M \colon \, \beta (x,\delta_{(S,\gamma)}) \geq 0 \right\}
	\end{align*} is a non-empty closed convex cone in $Y$. Moreover,
	\begin{displaymath}
	T \, \defeq \, \{ \omega \in \mathbb{R}[N] \mid
	\forall f \in N \colon \, \omega (f) \leq 0 \}
	\end{displaymath} is a non-empty convex cone in $Z$, which is closed
	in $Z$ by Proposition~\ref{proposition:bourbaki}(4).

	We now claim that statement~(1) of Corollary~\ref{corollary:farkas}
	is satisfied. Otherwise, there would exist some $x \in X$ such that
	$A(x) \in
	S$, $B(x) \in T$, and $\phi(x) < 0$, i.e., there would be some
	$\omega \in \mathbb{R}\!\left[\mathcal{O}_{D}^{(k)}\right]$ with
	$\sum_{f \in \supp(\omega)} \omega (f) = 0$ such that \begin{enumerate}
		\item[(a)] $\forall (S,\gamma) \in M \colon \, \sum_{f \in \supp
			(\omega)} \omega (f) \gamma (f(S)) \geq 0$,
		\item[(b)] $\forall f \in \mathcal{O}_{D}^{(k)} \setminus
		\mathcal{J}_{D}^{(k)} \colon \, \omega (f) \leq 0$,
		\item[(c)] $\sum_{f \in \supp (\omega)} \omega (f) \rho
		(f(x^{1},\ldots,x^{k})) < 0$,
	\end{enumerate} which would entail that $- \omega \in
	\mathrm{wPol}(\Gamma)$ (by~(a) and~(b)), while $- \omega$ would not
	weight-improve $\rho$ (due to~(c)), thus providing a contradiction.
	Therefore, item~(1) of Corollary~\ref{corollary:farkas} is satisfied,
	whence~(2) of the same holds, too.
	
	In order to show that $\rho \in \mathrm{\ell Expr}(\Gamma)$, let
	$\varepsilon > 0 $ and let $\mathcal F$ be a finite subset of
	$\mathcal{O}_{D}^{(k)}$. Without loss of generality, we may and will
	assume that $\mathcal{J}_{D}^{(k)}
	\subseteq \mathcal F$. By statement (2) of
	Corollary~\ref{corollary:farkas}, there exist $\mu \in S^{+}$ and $\nu
	\in T^{+}$ such that
	\begin{equation}\label{key}
	\forall f \in \mathcal{F} \ \forall e \in \mathcal{J}_{D}^{(k)}
	\colon \ \ \vert \phi(\delta_{f} - \delta_{e}) - (\mu(A(\delta_{f} -
	\delta_{e})) + \nu (B(\delta_{f} - \delta_{e}))) \vert \leq
	\frac{\varepsilon}{2} .
	\end{equation} Note that $\phi(\delta_{f} - \delta_{e}) = \rho
	(f(x^{1},\ldots,x^{k})) - \rho (e(x^{1},\ldots,x^{k}))$ whenever $f \in
	\mathcal
	F$ and $e \in \mathcal{J}_{D}^{(k)}$. Furthermore, since
	$(\mathbb{R}^{M},\sigma(\mathbb{R}^{M},\mathbb{R}[M]))^{\ast} \cong
	\mathbb{R}[M]$ (see Section~\ref{section:constructions}), we find
	$\lambda \in \mathbb{R}[M]$ such that \begin{displaymath}
	\forall x \in R^{M} \colon \quad \mu (x) \, = \,
	\sum\nolimits_{(S,\gamma) \in \supp (\lambda)} \lambda
	(S,\gamma)x(S,\gamma) .
	\end{displaymath} As $\mu \in S^{+}$, it follows that
	$\lambda(S,\gamma) = \mu (\delta_{(S,\gamma)}) \geq 0$ for every
	$(S,\gamma) \in M$. Furthermore, \begin{align*}
	\mu(A(\delta_{f} - \delta_{e})) \, = \, &\sum\nolimits_{(S,\gamma)
		\in \supp (\lambda)} \lambda (S,\gamma) (\gamma (f(S)) - \gamma(e(S))) \\
	= \,&\sum\nolimits_{(S,\gamma) \in \supp (\lambda)}
	\lambda(S,\gamma)\gamma (f(S)) - \sum\nolimits_{(S,\gamma) \in \supp
		(\lambda)} \lambda(S,\gamma)\gamma (e(S))
	\end{align*} for all $f \in \mathcal F$ and $e \in
	\mathcal{J}_{D}^{(k)}$. Moreover, $B(\delta_{f} - \delta_{e}) =
	\delta_{f}$ for all $f \in \mathcal{F}\setminus \mathcal{J}_{D}^{(k)}$
	and $e \in \mathcal{J}_{D}^{(k)}$. Now, for each $f\in \mathcal F$, let
	us define \begin{displaymath}
	C_{f} \defeq \begin{cases}
	\rho (f(x^{1},\ldots,x^{k})) - \sum\nolimits_{(S,\gamma)
		\in \supp (\lambda)} \lambda(S,\gamma)\gamma (f(S)) - \nu (\delta_{f}) &
	\text{if } f \notin \mathcal{J}_{D}^{(k)}, \\
	\rho (f(x^{1},\ldots,x^{k})) - \sum\nolimits_{(S,\gamma)
		\in \supp (\lambda)} \lambda(S,\gamma)\gamma (f(S)) & \text{if } f \in
	\mathcal{J}_{D}^{(k)};
	\end{cases}
	\end{displaymath} and for each $e \in \mathcal{J}_{D}^{(k)}$, let us
	define \begin{displaymath}
	K_{e} \, \defeq \, \rho (e(x^{1},\ldots,x^{k})) -
	\sum\nolimits_{(S,\gamma) \in \supp (\lambda)} \lambda(S,\gamma)\gamma
	(e(S)) .
	\end{displaymath} Evidently, $C_{e} = K_{e}$ for each $e \in
	\mathcal{J}_{D}^{(k)}$. Furthermore, from~\eqref{key} we infer that
	\begin{align*}
	\vert C_{f} - C_{e} \vert \, &= \, \left\lvert \phi(\delta_{f} -
	\delta_{e}) - (\mu(A(\delta_{f} - \delta_{e})) + \nu (\delta_{f}))
	\right\rvert \\
	& = \, \left\lvert \phi(\delta_{f} - \delta_{e}) - (\mu(A(\delta_{f}
	- \delta_{e})) + \nu (B(\delta_{f} - \delta_{e}))) \right\rvert  \, \leq
	\, \frac{\varepsilon}{2}
	\end{align*} for all $f \in \mathcal{F}\setminus \mathcal{J}_{D}^{(k)}$
	and $e \in \mathcal{J}_{D}^{(k)}$. Thus, \begin{equation}\label{better}
	\forall f \in \mathcal{F} \ \forall e \in \mathcal{J}_{D}^{(k)}
	\colon \qquad \left\vert C_{f} - K_{e} \right\vert \, \leq \,
	\frac{\varepsilon}{2}.
	\end{equation} Let $C \defeq \min_{f \in \mathcal F} C_{f}$. It
	follows from~\eqref{better} that $\vert C - K_{e} \vert \leq
	\frac{\varepsilon}{2}$ for each $e \in \mathcal{J}_{D}^{(k)}$,  which,
	by another application of~\eqref{better}, entails that $\vert C - C_{f}
	\vert \leq \varepsilon$ for every $f \in \mathcal F$. Let us prove that
	\begin{equation}\label{goal1}
	\forall f \in \mathcal{F} \colon \ \ \rho(f(x^{1},\ldots,x^{k})) \,
	\leq \, \sum\nolimits_{(S,\gamma) \in \supp (\lambda)} \lambda(S,\gamma)
	\gamma (f(S)) + C + \varepsilon.
	\end{equation} Clearly, for every $f \in \mathcal{J}_{D}^{(k)}$,
	\begin{align*}
	\sum\nolimits_{(S,\gamma) \in \supp (\lambda)} &\lambda(S,\gamma)
	\gamma (f(S)) + C \\
	& = \, \rho (f(x^{1},\ldots,x^{k})) - C_{f} + C \, \geq \, \rho
	(f(x^{1},\ldots,x^{k})) - \varepsilon .
	\end{align*} Moreover, for every $f \in N$, since $-\delta_{f} \in T$
	and $\nu \in T^{+}$, it follows that \begin{displaymath}
	\nu (\delta_{f}) \, = \, - \nu (-\delta_{f}) \, \leq \, 0 .
	\end{displaymath} Consequently, for each $f \in \mathcal F\setminus \mathcal{J}_{D}^{(k)}$,
	\begin{align*}
	\sum\nolimits_{(S,\gamma) \in \supp (\lambda)} &\lambda(S,\gamma)
	\gamma (f(S)) + C \\
	& \geq \, \sum\nolimits_{(S,\gamma) \in \supp (\lambda)}
	\lambda(S,\gamma) \gamma (f(S)) + \nu (\delta_{f}) + C  \\
	& = \, \rho (f(x^{1},\ldots,x^{k})) - C_{f} + C \, \geq \, \rho
	(f(x^{1},\ldots,x^{k})) - \varepsilon .
	\end{align*} This proves~\eqref{goal1}. Finally, we claim that
	\begin{equation}\label{goal2}
	\forall i \in \{ 1,\ldots,k \} \colon \ \left\lvert \rho(x^{i}) -
	\left( \sum\nolimits_{(S,\gamma) \in \supp (\lambda)} \lambda(S,\gamma)
	\gamma \!\left( e_{i}^{(k)}(S) \right) + C \right) \right\rvert \leq
	\varepsilon.
	\end{equation} Indeed, if $i \in \{ 1,\ldots,k \}$, then we consider
	$e \defeq e_{i}^{(k)}$ and observe that \begin{align*}
	\sum\nolimits_{(S,\gamma) \in \supp (\lambda)} \lambda(S,\gamma)
	\gamma (e(S)) + C \, & \leq \, \sum\nolimits_{(S,\gamma) \in \supp
		(\lambda)} \lambda(S,\gamma) \gamma (e(S)) + K_{e} +
	\frac{\varepsilon}{2} \\
	& = \, \rho (e(x^{1},\ldots,x^{k})) + \frac{\varepsilon}{2} \, =
	\, \rho (x^{i}) + \frac{\varepsilon}{2} .
	\end{align*} Combined with~\eqref{goal1} above and the fact that $e
	\in \mathcal F$, this implies~\eqref{goal2}. This proves that $\rho \in
	\mathrm{\ell Expr}(\Gamma)$, as desired. \end{proof}

\section*{Conclusion}

We have given an algebraic characterisation of the local expressive power 
of finite-valued constraint languages over domains of arbitrary cardinality. 

In the finite-domain setting, the tractability of the VCSP for a valued structure implies the tractability of the VCSP for its expressive power, and the NP-hardness of the VCSP for a valued structure can be inferred from the expressibility of cost functions that encode NP-hard problems, such 
as \textsc{Max-Cut} and 3-SAT. This result can be extended to the infinite-domain finite-valued setting under some additional assumptions, as we have seen.
Because of this complexity reduction, in the finite-domain setting the characterisation of the expressive power in terms of weighted polymorphisms implies a characterisation of the computational complexity in terms of 
algebraic properties of the valued structure.

It is interesting to understand whether such an algebraic characterisation of the computational complexity of VCSPs can be extended to the infinite-domain case. We know that for general-valued CSPs weighted polymorphisms do not capture the computational complexity, since already the complexity of an infinite-domain $\{0,+\infty\}$-valued constraint language is not captured by its polymorphisms~\cite[Section~4.1]{Bodirsky2017ConstraintSP}. However, it makes sense to question the 
relationship between the expressive power and the set of weighted polymorphisms for finite-valued CSPs over infinite domains. 

The present paper represents the first progress to answer this question; in view of our result, investigating the relationship between the expressive power and the set of weighted polymorphisms reduces to understanding 
whether the containment of the expressive power into the local expressive 
power of a finite-valued structure is strict.

\section*{Acknowledgements}
The authors are grateful to Manuel Bodirsky for suggesting the topic of the present paper and for the many helpful discussions. Furthermore, the authors would like to express their sincere gratitude towards the anonymous referee for numerous comments and suggestions on an earlier version of this manuscript that helped to significantly improve the presentation.\;

Caterina Viola has been supported by the Deutsche Forschungsgemeinschaf (DFG) Graduiertenkolleg 1763 (QuantLA) and has received funding from the European Research Council (ERC) under the European Union's Horizon 2020 research and innovation programme (grant agreement No 681988, CSP-Infinity).



%%       \citet{<label>} ==> Jones et al. [21]
%%       \citep{<label>} ==> [21]
%%



\bibliographystyle{unsrtnat}
\bibliography{locexppwr}

\end{document}